\documentclass[a4paper,12pt]{article}

\usepackage{amssymb}

\usepackage{amsmath}
\usepackage{amsfonts}
\usepackage{amscd}
\usepackage[english]{babel}

\title{Cayley structures on $S^6$ as a twistor bundle sections}
\author{Natalia Daurtseva\\ Russia, Kemerovo State University\footnote{natali0112@ngs.ru}}
\date{ }
\begin{document}
\maketitle

The nearly K\"{a}hler structures on the 6-sphere, as a twistor
bundle sections are researched. We show that for any point of
twistor bundle there exists an 1-parametric family of sections,
passing through the point, which give nearly K\"{a}hler structures
on the round sphere. Some properties of those sections are found.


\section*{Introduction} % Раздел без номера

We will talk about Cayley structures on the 6-sphere. With round
metric they define nearly K\"{a}hler structures on $S^6$.

{\bf Definition 1.1.}\textit{Nearly K\"{a}hler manifold is an
almost Hermitian manifold $(M,$ $g,$ $J$, $\omega)$ with the
property that $(\nabla_XJ)X=0$ for all tangent vectors $X$, where
$\nabla$ denotes the Levi-Civita connection of $g$. If
$\nabla_XJ\neq 0$ for any non-zero vector field $X$,
$(M,g,J,\omega)$ is called strictly nearly K\"{a}hler.}

Nearly K\"{a}hler geometry comes from the concept of weak holonomy
introduced by A. Gray in 1971 [1], this geometry corresponds to
weak holonomy $U(n)$. Lie group $U(n)$ is a structure group of
almost Hermitian manifold, if holonomy group is equal to $U(n)$
too, then this manifold is K\"{a}hler. In weak case  \cite{G1}
almost Hermitian manifold with weak holonomy group $U(n)$ is
nearly K\"{a}hler. The class of nearly K\"{a}hler manifolds
appears naturally as one of the sixteen classes of almost
Hermitian manifolds described by the Gray-Hervella classification
[2].

At 2002 Nagy P.-A. has proved \cite{Nagy} that every compact
simply connected nearly K\"{a}hler manifold $M$ is isometric to a
Riemannian product $M_1\times\dots M_k$, such that for each $i$,
$M_i$ is a nearly K\"{a}hler manifold belonging to the following
list: K\"{a}hler manifolds, naturally reductive 3-symmetric
spaces, twistor spaces over compact quaternion-K\"{a}hler
manifolds with positive scalar curvature, and nearly K\"{a}hler
6-manifold. This is one of the reasons of interest to the nearly
K\"{a}hler 6-manifolds.

For 6-dimensional nearly K\"{a}hler Riemannian homogeneous
manifolds we have the following classification:

{\bf Theorem (J.P.Butruille \cite{But})} {\it Nearly K\"{a}hler,
6-dimensional, Riemannian homogeneous manifolds are isomorphic to
a finite quotient of $G/H$ where the groups $G$, $H$ are given in
the list:
\begin{itemize}
\item[--] $G=SU(2) \times SU(2)$ and $H=\{1\}$
\item[--] $G=Sp(2)$ and $H=SU(2)U(1)$. Then, $G/H \simeq \mathbb{C}P^3$, the 3-dimensional complex projective space.
\item[--] $G=SU(3)$, $H=U(1) \times U(1)$ and $G/H$ is the space of flags of $\mathbb{F}^3$.
\item[--] $G=G_2$ and $H=SU(3)$. In this case $G/H$ is the round 6-sphere.
\end{itemize}
}

And we don't know an example of the nonhomogeneous nearly
K\"{a}hler manifold.

For the first three manifolds from the list above there exists
only one nearly K\"{a}hler structure. But in case of the round
6-sphere one has an infinite number of those structures.

{\bf Theorem (M. Verbitsky \cite{Verb})} {\it Let $(M^6, I, g)$ be
a nearly K\"{a}hler manifold. Then the almost complex structure is
uniquely determined by the Riemannian structure, unless M is
locally isometric to a 6-sphere.}

All almost complex structures, which gives the nearly K\"{a}hler
structures together with the round metric are Cayley structures on
$S^6$.

\section{Cayley structures on $S^6$}

Let's look $\mathbb{R}^7$ as the space of purely imaginary Cayley
numbers with inner product
$$
\langle x,y\rangle=-\mbox{Re}(xy),\qquad\forall x,y\in\mathbb{R}^7
$$
Then, sphere $S^6$ is a subset of unit length numbers in
$\mathbb{R}^7$. Inner product $<,>$ induces round metric $g_0$ on
$S^6$.

Let $R_x:\mathbb{O}\longrightarrow\mathbb{O}$ denotes the right
Cayley multiplication by $x$, $R_xy=y\cdot x$. The right
multiplication defines standard Cayley structure $\mathcal{J}$ on
the $S^6$:
$$
\mathcal{J}=\{(x,R_x):x\in S^6\}
$$
For any vector $y\in T_xS^6$, $\mathcal{J}_x(y)=R_x(y)=y\cdot x$,
$\forall x\in S^6$. By the multiplication properties \cite{CG},
one can show that $\mathcal{J}_x(y)\in T_xS^6$ and
$$
\mathcal{J}_x^2(y)=R_x(R_x(y))=(y\cdot x)\cdot x=y\cdot(x\cdot
x)=-y
$$
So, the Cayley structure is an almost complex structure on $S^6$.

Another Cayley structures are defined from the following way. Let
$A\in O(7)$ is any orthogonal transformation, then \cite{CG} one
can define structure $\mathcal{J}^A$:
$$\mathcal{J}^A=\{(x,A^{-1}R_{A(x)}A):x\in S^6\},$$
Obviously,
$\mathcal{J}^A_x(y)=A^{-1}R_{A(x)}A(y)=A^{-1}(A(y)\cdot A(x))$,
$\forall x\in S^6$, $\forall y\in T_xS^6$.

Structure $\mathcal{J}^A$ is equal to standard one if and only if
$A$ is automorphism of Cayley algebra. So, the space of all Cayley
structures is homogeneous space $O(7)/G_2$.

{\bf Remark 1.} {\it Further we will speak about almost complex
structures which define the standard orientation, induced from
$\mathbb{R}^7$, only. In this case the space of those structures
is $SO(7)/G_2\cong\mathbb{R}P^7$.}

All these structures are nearly K\"{a}hler on the homogeneous
Riemannian manifold $S^6=G_2/SU(3)$. Every embedding $G_2$ in
$SO(7)$ defines on the sphere some nearly K\"{a}hler structure
associated with round metric.

\section{Twistor bundle over $S^6$}

Let $(Z^+(S^6,g_0),\pi)$ is a twistor bundle over $S^6$:
$$Z^+(S^6,g_0)=\{(x,J_x):x\in S^6, J_x:T_xS^6\longrightarrow T_xS^6, J_x^2=-1$$
$$\mbox{and define the standard orientation, induced from }
\mathbb{R}^7\}$$
$$\pi:Z^+(S^6,g_0)\longrightarrow S^6,\ \pi(x,J_x)=x$$

Each fiber of bundle is 6-dimensional and diffeomorphic to
$SO(6)/U(3)\cong\mathbb{C}P^3$. Thus, the space $Z^+(S^6,g_0)$ is
12-dimensional. Any almost complex structure $\mathcal{J}$ on
$(S^6,g_0)$ is a smooth section of the bundle and
$\mathcal{J}(S^6)$ is smooth submanifold in $Z^+(S^6,g_0)$.

\vspace{3mm}

{\bf Lemma 1.}\label{L1} {\it For any point $(p, J_p)\in
Z^+(S^6,g_0)$ there exists Cayley structure $\mathcal{J}^A$, with
$\mathcal{J}_p^A=J_p$.}

{\bf Proof.} Let $p\in S^6$ is any point on the sphere, and
$J_p\in\pi^{-1}(p)$ is arbitrary. Let $\mathcal{J}_p$ is value of
standard Cayley structure at $p$. Then one can find corresponding
orthogonal transformation $A\in SO(6)$ of $T_pS^6$ with
$\mathcal{J}_p=AJ_pA^{-1}$. Look $A$, as rotation around the point
$p$ in $\mathbb{R}^7$. Then $A$ defines Cayley strucutre
$$
\mathcal{J}^A=\{(x,A^{-1}\mathcal{J}_{A(x)}A):x\in S^6\}
$$
At point $p$ we have
$\mathcal{J}^A_p(y)=A^{-1}\mathcal{J}_{A(p)}A(y)=A^{-1}\mathcal{J}_pA(y)=J_p(y)$
for any vector $y\in T_pS^6$. Lemma is proved.

\vspace{3mm}

Dimension of the fiber is equal to 6, but dimension of the Cayley
structures space $SO(7)/G_2$ is equal to 7. So it would be natural
to suggest that in any fiber we have points where more then one
Cayley structure pass through.

\vspace{3mm}

{\bf Theorem 1.}\label{T1} {\it For any point $(p, J_p)\in
Z^+(S^6,g_0)$  there exists an 1-parametric family of twistor
bundle sections, passing through the point, which give nearly
K\"{a}hler structures on $(S^6,g_0)$.}

\vspace{3mm}

{\bf Proof.} Let $(p,J_p)\in Z^+(S^6,g_0)$ is any point of twistor
bundle over $S^6$. Let $\mathcal{J}^A$ is corresponding to $J_p$
(as in Lemma 1), $\mathcal{J}^A_p=J_p$, $A\in SO(6)$. Group
$SO(6)$ in this case is a subgroup in $SO(7)$ of rotations around
the $p$. Let $\lambda=(\cos\varphi+i\sin\varphi)E\in U(3)$,
$\varphi\in [0;2\pi/3)$ is unitary transformation in the plane
$T_pS^6$, commuting with $\mathcal{J}_p$. This transformation
$\lambda\notin SU(3)=G_2\cap SO(6)$, while $\varphi\neq 0$, then
$$
\mathcal{J}^{\lambda A}_p=A^{-1}\lambda^{-1}\mathcal{J}_p\lambda
A=A^{-1}\mathcal{J}_pA=\mathcal{J}^A_p=J_p,
$$
but $\mathcal{J}^{\lambda_1 A}\neq\mathcal{J}^{\lambda_2 A}$ for
any different $\varphi_1$ and $\varphi_2$, as in this case
$\lambda_1\lambda_2^{-1}\notin G_2$. Theorem is proved.

\vspace{3mm}

{\bf Remark 2. }{\it Denote 1-parametric family of Cayley
structures through the point $(p,J_p)\in Z^+(S^6,g_0)$ as
$\{p,J_p\}$.}

\vspace{3mm}

{\bf Lemma 2.}{\it Let $(p, J_p)\in Z^+(S^6,g_0)$ then, any two
Cayley structures $\mathcal{J}^A$ and $\mathcal{J}^B\in\{p,J_p\}$
have exactly two intersection points.}

\vspace{3mm}

{\bf Proof.} Let $\mathcal{J}^A$ and $\mathcal{J}^B\in\{p,J_p\}$.
By lemma 2 and proof of theorem 1 one can think about
transformations $A$ and $B$ as a rotations in $\mathbb{R}^7$ with
respect to $p$. Moreover $B=\lambda A$, where
$\lambda=(\cos\varphi+i\sin\varphi)E\in U(3)$, $\varphi\in
(0;2\pi/3)$. Rotation around the point $p$ in $\mathbb{R}^7$ keeps
point $p$ and $-p$ invariant. Obviously, that values of structures
$\mathcal{J}^A$ and $\mathcal{J}^B$ at $-p$ are equal.

Suppose that one can find a point $q\neq\pm p$ on the sphere, and
point $J_q\in\pi^{-1}(q)$, with
$\mathcal{J}^A_q=\mathcal{J}^{A\lambda}_q=J_q$. Then by definition
of Cayley structure and  proof of theorem 1:
$$
(\lambda A)^{-1}R_{\lambda A(q)}(\lambda A)(y)= A^{-1}R_{A(q)}A(y)
$$
for all $y\in T_qS^6$.
$$
\lambda^{-1}R_{\lambda q'}(\lambda y')=R_{q'}(y')
$$
where $q'=Aq$, $y'=Ay$. As $q'\neq\pm p,0$ and $y'$ is arbitrary
tangent vector at point $q'$, then the equality is possible just
if $\lambda\in G_2$, this contradicts with choice of $\lambda$.

\vspace{3mm}

Let $(p,J_p)$ is some point of the space $Z^+(S^6,g_0)$, and
$\mathcal{J}\in\{p,J_p\}$ is any section of the bundle, through
the point. Then tangent space to $\mathcal{J}(S^6)$ is:
$$
T_{(p,J_p)}\mathcal{J}(S^6)=\{(X,K)\in T_pS^6\times
T_{J_p}\pi^{-1}(p): X=\frac{dx_t}{dt}|_{t=0},\
K=\frac{dJ_{x_t}}{dt}|_{t=0};\ x_t:I\longrightarrow S^6,\ x_0=p\}
$$

\vspace{3mm}

{\bf Lemma 3. }{\it For any point $(p,J_p)\in Z^+(S^6,g_0)$ any
two sections $\mathcal{J}^A$ and $\mathcal{J}^{\lambda
A}\in\{p,J_p\}$ intersect transversely.}

\vspace{3mm}

{\bf Proof.} Tangent spaces $T_{(p,J_p)}\mathcal{J}^A(S^6)$ and
$T_{(p,J_p)}\mathcal{J}^{\lambda}(S^6)$ are:
$$
T_{(p,J_p)}\mathcal{J}^A(S^6)=\{(X,K):X\in T_pS^6;
K=A^{-1}R_{AX}A\}
$$
$$
T_{(p,J_p)}\mathcal{J}^{A\lambda}(S^6)=\{(X,K_{\lambda}):X\in
T_pS^6; K_{\lambda}=A^{-1}\lambda^{-1}R_{\lambda A X}\lambda A\}
$$
$(X,K)=(X,K_{\lambda})$ if and only if
$$
K=K_{\lambda}
$$
$$
A^{-1}R_{AX}A(y)=A^{-1}\lambda^{-1}R_{\lambda A X}\lambda A(y)
$$
for any vector $y\in T_p(S^6)$
$$
R_{X'}(y')=\lambda^{-1}R_{\lambda X'}\lambda (y'),
$$
where $X'=AX$, $X'\neq\pm p; 0$, $y'=A(y)$. This is possible just
in case of $\lambda\in G_2$, but this contradicts to choice of
$\lambda$.

\vspace{3mm}

{\bf Lemma 4. }{\it Any Cayley structure $\mathcal{J}^A$, passing
through a point $(q,J_q)\in Z^+(S^6,g_0)$ is in family
$\{q,J_q\}$.}

\vspace{3mm}

{\bf Proof.} Let Cayley structure $\mathcal{J}^A$ comes through
the point $(q,J_q)$, i.e. $\mathcal{J}^A_q=J_q$. Then by theorem 1
and lemma 1 $\mathcal{J}^A\in\{q,J_q\}$ just in case of existence
the transformation $G\in G_2$, with $GA(q)=q$. Lie group $G_2$
acts transitively on sphere, so for any points $q, A(q)\in S^6$
one can find transformation $G\in G_2$, for which $GA(q)=q$.

\vspace{3mm}

{\bf Corollary 1.} {\it If any two Cayley structures are
$\mathcal{J}^A$ and $\mathcal{J}^B$ intersect at a point
$(q,J_q)\in Z^+(S^6,g_0)$, then $\mathcal{J}^A,
\mathcal{J}^B\in\{q,J_q\}$.}

\section{Nearly K\"{a}hler $SU(3)$-structures on $S^6$}

We can see the above 1-parametric families using the differential
forms point of view. Two forms $(\omega,\psi)$ on $M^6$, where
$\omega$ - is nondegenerate skewsymmetric  2-form, with stabilizer
$Sp(3,\mathbb{R})$, and $\psi$ - 3-form with stabilizer
$SL(3,\mathbb{C})$ in each point, with some additional conditions
define $SU(3)$-structure $M$. Really, 3-form $\psi$ gives
reduction to $SL(3,\mathbb{C})$, and $\omega$ to
$Sp(3,\mathbb{R})$. Group $SU(3)$ is intersection
$Sp(3,\mathbb{R})\cap SL(3,\mathbb{C})$, but to get
$SU(3)$-structure we need in two more algebraic conditions. To
start remind that 3-form on $M$ with $SL(3,\mathbb{C})$ as
stabilizer defines almost complex structure $J_{\psi}$ by formulas
\cite{Hitchin}:
$$
K(X)=A(\iota_X\psi\wedge\psi)
$$
where $A:\Lambda^5\longrightarrow TM\otimes\Lambda^6$ is
isomorphism, induced by exterior product
($\iota_{A(\varphi)}Vol=\varphi$). Then,
$\tau(\psi)=\frac16\mbox{tr}K^2$, is section
$\Lambda^6\otimes\Lambda^6$ and $K^2=Id\otimes\tau(\psi)$. It is
known, that $SL(3,\mathbb{C})$ is stabilizer of 3-form $\psi$ at
each point, if and only if $\tau(\psi)<0$. In this case, form
$\psi$ defines almost complex structure on $M$
$$
J_{\psi}=\frac{1}{\kappa}K
$$
where $\kappa=\sqrt{-\tau(\psi)}$.

The first of above algebraic properties is:

$$\omega\wedge\psi=0$$

This is the property of $\omega$ to be of type (1,1) with respect
to $J$. Second condition is  positivity of form $\omega(X,JY)>0$.
In case, if the first and the second properties are hold, then
forms $(\omega,\psi)$ define $SU(3)$ structure, and give almost
complex structure $J_{\psi}$, metric $g(X,Y)=\omega(X,J_{\psi}Y)$,
such that $g(J_{\psi}X,J_{\psi}Y)=g(X,Y)$. In this case almost
Hermitian structure $(g,J_{\psi})$ is nearly K\"{a}hler if the
following conditions are hold \cite{RC}:
$$\left\{\begin{array}{l} \psi=3d\omega\\
d\phi=-2\mu\omega\wedge\omega
\end{array}\right.,$$ where $\iota_{JX}\phi=\iota_X\psi$.

Now, the cone of $(M,g)$ is the Riemannian manifold $(\overline M,
\overline g)$ where $\overline M = M \times \mathbb{R}^+$ and
$\overline g = r^2 g + dr^2$ in the coordinates $(x,r)$. Let
define a section $\rho$ of $\Lambda^3 \overline M$ by
$$
\rho=r^2dr\wedge\omega+r^3\psi,
$$
$\rho$ is a generic 3-form, inducing a $G_2$-structure on the
7-manifold $\overline M$ such that $\overline g$ is the metric
determined by $\rho$, given the inclusion of $G_2$ in $SO(7)$.
Moreover $\nabla^{\overline g} \rho = 0$, where $\nabla^{\overline
g}$ is the Levi-Civita connection of $\overline g$. In other
words, the holonomy of $(\overline M,\overline g)$ is contained in
$G_2$.

Reciprocally, a parallel, generic 3-form on $\overline M$ can
always be written $\rho=r^2dr\wedge\omega+r^3\psi$, where
$(\omega,\psi)$ define a nearly K\"{a}hler $SU(3)$-structure on
$M$.

The Riemannian cone of the 6-sphere is the Euclidean space
$\mathbb{R}^7$. According to what precedes, nearly K\"{a}hler
structure on $S^6$, compatible with $g_0$, define a parallel or
equivalently, a constant 3-form on $\mathbb{R}^7$.

Let $x\in S^6$, and $(\omega_x,\psi_x)$ is some $SU(3)$ structure
on $T_xS^6$. This structure is used to construct nearly K\"{a}hler
structure on the sphere. Define constant 3-form on $\mathbb{R}^7$:
$$
\rho_{(x,1)}=dr\wedge\omega_x+\psi_x.
$$
Then, $\rho$ is parallel for the Levi-Civita connection of the
flat metric and in return, $\omega=\iota_{\partial/\partial
r}\rho$ and $\psi=\frac13 d\omega$ determine a nearly Kдhler
structure on $S^6$ whose values at $x$ coincide with $\omega_x$,
$\psi_x$, consistent with our notations \cite{But}.

Let $\Psi_x=dz^1\wedge dz^2\wedge dz^3$ is complex volume form on
$T_xS^6$, then $\psi_x=\mbox{Re}\Psi_x$ defines standard complex
structure $J_{0_x}$ on $T_xS^6=\mathbb{R}^6$. Let
$\lambda=\cos\varphi+i\sin\varphi$,
$\varphi\in[0;\frac{2\pi}{3})$, consider family of forms
$\Psi_{\lambda_x}=\lambda^3\Psi_x$. Identify $\lambda$ with
unitary transformation $\lambda\cdot Id\in U(3)$, embedded into
$SO(6)$ by standard way. Then
$\psi_{\lambda_x}(X,Y,Z)=\psi_x(\lambda X,\lambda Y,\lambda Z)$
and $K_{\lambda}(\lambda X)=\lambda K(X)$, $J_{\lambda_x}=\lambda
J_x\lambda^{-1}=J_x$. So, 1-parametric family of
$SU(3)$-structures $(\omega_x, \mbox{Re}(\lambda \Psi_x))$ defines
family of nearly K\"{a}hler structures $J_{\lambda}$, equal to
each other at point $x$.

Any another $SU(3)$ structure on $T_xS^6$ is defined by couple of
forms $(\omega_A,\psi_A)$, where $A\in SO(6)$,
$\omega_A(X,Y)=\omega(A^{-1}X,A^{-1}Y)$,
$\psi_A(X,Y,Z)=\psi(AX,AY,AZ)$. For those forms $K_A(AX)=AK(X)$
and $J_{A_x}=AJ_xA^{-1}$. So, at any $x\in S^6$ $SU(3)$ structures
$(\omega_x,\psi_x)$ defines all possible almost complex structures
on $T_xS^6$, for any $J_x$ there exists 1-parametric family of
nearly K\"{a}hler structures $J_{\lambda}$ on $S^6$, with
$J_{\lambda_x}=J_x$.

{\bf Remark 3} {\it If $p$ is some fixed point on the 6-sphere,
then any $SU(3)$ structure on $T_xS^6$ defines nearly K\"{a}hler
structure on $S^6$. Set of those $SU(3)$ structures on
6-dimensional vector space is 7-dimensional homogeneous space
$SO(6)/SU(3)\cong SO(7)/G_2\cong\mathbb{R}P^7$. The fiber of
twistor bundleis diffeomorphic to $SO(6)/U(3)$. Through
"difference" between $U(3)$ and $SU(3)$ the above families of
nearly K\"{a}hler structures are arise.}

\end{document}